# The Kruskal-Katona Theorem for Graphs

by Robert Cowen,
Queens College, CUNY

Abstract. In graph theory, knowing the number of complete subgraphs with $r$ vertices that a graph $g$ has, limits the number of its complete subgraphs with $s$ vertices, for $s > r$. A useful upper bound is provided by the Kruskal-Katona theorem, but this bound is often not tight. In this note, we add to the known cases where this bound is tight and also investigate cases where it is not. Finally we look at some useful techniques for actually finding the numbers of complete subgraphs of a graph.

1. Introduction. In graph theory, knowing the number of complete subgraphs with $r$ vertices that a graph $g$ has, limits the number of its complete subgraphs with $s$ vertices, for $s > r$. A useful upper bound was given by J. Kruskal in [5], rediscovered by G. Katona in [4]. This, so-called Kruskal-Katona bound is not always tight(see [2]); however, there are obvious cases when it is a tight bound, for example, complete graphs. In [1], B. Bollobás proved a theorem that provided some additional graphs where these bounds are tight.

Definition. If $x$ is a non-negative integer its $r$-canonical representation is $\binom{n}{r} + \binom{m}{r-1} + \ldots + \binom{u}{r-j}$, where $n$ is chosen as large as possible with $\binom{n}{r} < x$, then $m$ is chosen as large as possible with $\binom{m}{r-1} < (x - \binom{n}{r})$, ..., etc., and, $x = \binom{n}{r} + \binom{m}{r-1} + \ldots + \binom{u}{r-j}$. We shall denote the r-canonical representation of $x$ by $[x]^r$. If $r < s$, let $[x]^r_s$ be the result of replacing $r$ by $s$ in the $r$-canonical representation of $x$, that is, $[x]^r_s = \binom{n}{s} + \binom{m}{s-1} + \ldots + \binom{u}{s-j}$.

It is easy to see that $n > m > \ldots > u$. Also, any binomial coefficient whose top entry is less than its bottom entry is taken to be 0.

Definition. $K_r$ shall denote the complete graph with $r$ vertices, that is the graph with $r$ vertices with an edge between any pair of its vertices. If $g$ is a graph and $r > 1$, $k_r(g)$ will denote the number of its $K_r$ subgraphs.

The following theorem is a consequence of Theorem 1 in Kruskal[5].

Theorem 1. If $r < s$ and $g$ is any graph with $k_r(g) \leq x$, then $k_s(g) \leq [x]^r_s$.

The upper bound in Theorem 1 is known as the Kruskal-Katona bound and it is not always tight (see [2]).

Since this upper bound isn't always tight, we shall introduce the following notation due to B. Bollobás [1].

Definition. If $r < s$, $k_s(k_r \leq x)$ is the maximum number of $K_s$ subgraphs in a graph that has at most $x$, $K_r$ subgraphs.

It is easy to see that $h(x) = k_s(k_r \leq x)$ is a non-decreasing function of $x$.

The following then is just a restatement of Theorem 1.
Theorem 1. If $r < s$, then $k_s(k_r \leq x) \leq [x]_s^r$

2. Bollobás' theorem and extensions. It is known that equality is not always achieved in Theorem 1; for example, see [2]. However, the following result is due to B. Bollobás [ 1 ].

Theorem 2. Suppose $r < s < n$. If $[x]^r = \binom{n}{r} + \binom{m}{r-1}$, we have $k_s(k_r \leq x) = \binom{n}{s} + \binom{m}{s-1} = [x]_s^r$.

Our next theorem extends this result of Bollobás'.

Theorem 3. Suppose $r < s$ and $[x]^r = \binom{n}{r} + \binom{m}{r-1} + \binom{t}{r-2}$. If $\binom{t}{r-2} = \binom{w}{r-1}$, with $(s-2) > t$ and $(s-1) > w$, then $k_s(k_r \leq x) = [x]_s^r = \binom{n}{s} + \binom{m}{s-1}$. Moreover, $k_s(k_r \leq y) = [y]_s^r = \binom{n}{s} + \binom{m}{s-1}$ for all $y$ with $\binom{n}{r} + \binom{m}{r-1} \leq y \leq x$.

Proof. We show first that $k_s(k_r \leq x) = \binom{n}{s} + \binom{m}{s-1}$.

Since $r < s$, $k_s(k_r \leq x) \leq [x]_s^r$, by Theorem 1. However, $[x]_s^r = \binom{n}{s} + \binom{m}{s-1}$, because $\binom{t}{s-2} = 0$, since $(s-2) > t$. Hence $k_s(k_r \leq x) \leq \binom{n}{s} + \binom{m}{s-1}$.

We now construct a graph $g$ with $x = k_r(g)$ and $k_s(g) = [x]_s^r = \binom{n}{s} + \binom{m}{s-1}$. Let $g$ be the graph obtained by adding an external vertex A to the complete graph $K_n$ and then connecting A to $m$ vertices of $K_n$ and also adding an external vertex B and then connecting B to $w$ vertices of $K_n$. The number of $K_r$ subgraphs of $g$ is $\binom{n}{r}$ (those in $K_n$), plus $\binom{m}{r-1}$ (those containing the vertex A together with $(r-1)$ of the vertices that connect A to $K_n$), plus $\binom{t}{r-2} = \binom{w}{r-1}$ (those containing vertex B together with $(r-1)$ of the vertices connecting B to $K_n$. Adding up the three terms, gives $x$ $K_r$ subgraphs in $g$. We claim that $g$ has $\binom{n}{s} + \binom{m}{s-1}$ complete subgraphs, $K_s$. Surely $g$ has all the $K_s$ subgraphs of $K_n$ and, in addition, those $K_s$ subgraphs of $g$ consisting of the induced graphs whose vertex sets consist of vertex A, together with $(s-1)$ of the m vertices of $K_n$ connected to A; there is no contribution from the subgraphs containing B, since $(s-1) > w$. Hence $g$ has $\binom{n}{s} + \binom{m}{s-1}$ complete subgraphs, $K_s$, as claimed. There-fore, $k_s(k_r \leq x) \geq \binom{n}{s} + \binom{m}{s-1}$. Therefore, $k_s(k_r \leq x) = \binom{n}{s} + \binom{m}{s-1}$, as claimed.

Finally, $k_s(k_r \leq \binom{n}{r} + \binom{m}{r-1}) = \binom{n}{s} + \binom{m}{s-1}$, by the aforementioned result of Bollobás [1] and $k_s(k_r \leq u)$ is a non-decreasing function of $u$; hence, $k_s(k_r \leq y) = [y]_s^r = \binom{n}{s} + \binom{m}{s-1}$, for all $y$ with $\binom{n}{r} + \binom{m}{r-1} \leq y \leq x$.

<u>Corollary 1</u>. Let $u > 1$, $r = u + 1$ and suppose $s \geq 2u + 2$. If $[x]^r = \binom{n}{r} + \binom{m}{r-1} + \binom{2u-1}{r-2}$, then $k_s(k_r \leq y) = \binom{n}{s} + \binom{m}{s-1}$, for $\binom{n}{r} + \binom{m}{r-1} \leq y \leq x$.

<u>Proof</u>. Let $w = t = (2u - 1)$. Then $\binom{t}{r-2} = \binom{2u-1}{u-1} = \binom{2u-1}{u} = \binom{w}{r-1}$. Also, $s - 2 \geq 2u > t$ and $s - 1 > w$. Hence the result follows from Theorem 2.

The next result follows immediately from Corollary 1.

<u>Corollary 2</u>. Let $u > 1$, $r = u + 1$ and suppose $s \geq 2u + 2$. Then there is a sequence S of consecutive integers of length $\binom{2u-1}{r-2} + 1$ such that if $y \in S$, then $k_s(k_r \leq y) = [y]_s^r$, the Kruskal-Katona bound.

<u>Example</u>. Suppose $u = 4$. Then $r = 5$, and take $t = w = 7$. Also, if $s = 10$, then $s - 2 > t$ and $s - 1 > w$. Therefore, Corollary 1 of Theorem 2 implies that $k_{10}(k_5 \leq \binom{n}{5} + \binom{m}{4} + \binom{7}{3}) = \binom{n}{10} + \binom{m}{9}$, and we can now chose any $n, m$, with $n > m > 9$. For example, if we take $n = 11$ and $m = 10$, $\binom{11}{5} + \binom{10}{4} + \binom{7}{3} = 707$ and we get, $k_{10}(k_5 \leq 707) = \binom{11}{10} + \binom{10}{9} = 21$. The graph constructed with these numbers is shown in Figure 1.

3. <u>Removing Edges from Complete Graphs</u>. Suppose now that instead of adding to a complete graph a vertex and edges, we start with a compete graph and subtract edges from it. Suppose, for example, we subtract from the complete graph, $K_n$, a vertex and $p$ edges connected to that vertex, where $p < n$. In effect, we are left with a graph isomorphic to the complete graph, $K_{n-1}$ with a new vertex joined by m = n - p edges to $K_{n-1}$. Thus we get the following as a corollary to Theorem 2.

<u>Theorem 4</u>. If graph $g$ is obtained by deleting $p$ edges from the complete graph $K_n$, $p < n$, where all $p$ edges are incident to a single vertex, then $k_s(g) = \binom{n}{s} + \binom{n-p}{s-1}$.

If not all the edges to be deleted are incident with a single vertex, the situation is not as simple; for example, suppose just two edges are deleted, but these edges are not incident. These graphs, complete graphs with two non-incident edges removed, are the same as the Turan graphs, $T(n, n-2)$. (The Turan graph $T(n, k)$ is the graph formed by partitioning a set of $n$ vertices into $k$ subsets with sizes as equal as possible (differing by at most 1) and connecting two vertices by an edge if and only if they belong to different sets of the partition. This follows because the partition sets of $T(n, n-2)$ can be taken to be: {1, 2}, {3, 4}, {5}, {6}, …, {n}; so only edges (1, 2) and (3, 4) are forbidden.)

Next, consider the complete graph on 12 vertices with edges (1,2) and (3,4) deleted( see figure 2). Let's compute the numbers of $K_3$ and $K_4$ subgraphs of g.

There are 200 $K_3$ subgraphs ($\binom{12}{3} - 2(10)$) and 406 $K_4$ subgraphs ($\binom{12}{4} - 2\binom{10}{2} + 1$) ; however, the Kruskal-Katona bound for the $K_4$ subgraphs, for **any** graph with 200 $K_3$ subgraphs is $[200]_4^3 = 407$ which is only one more than the actual number of $K_4$ subgraphs in g! This is true in general, that is, if we consider, $K_n$ with two non incident edges deleted, there are ($\binom{n}{3} - 2(n-2)$) $K_3$ subgraphs and ($\binom{n}{4} - 2\binom{n-2}{2} + 1$) $K_4$ subgraphs. If $x = (\binom{n}{3} - 2(n-2))$, then it is not difficult to show that $x = \binom{n-1}{3} + \binom{n-4}{2} + \binom{n-5}{1}$ and hence, $[x]_4^3 = \binom{n-1}{4} + \binom{n-4}{3} + \binom{n-5}{2}$. Using a computer algebra program, such as Mathematica or using the Pascal identity judiciously, it can be shown that

$(\binom{n}{4} - 2\binom{n-2}{2} + 1) = \binom{n-1}{4} + \binom{n-4}{3} + \binom{n-5}{2} - 1$, n> 6; that is, the number of $K_4$ subgraphs in T(n, n - 2) is one less than the Kruskal-Katona bound, $[x]_4^3$, where x is the number of $K_3$ subgraphs of T(n, n - 2).

Thus, $k_4(k_3 \le x) \ge [x]_4^3 - 1$, when $x = (\binom{n}{3} - 2(n-2))$, n > 6. Since it is always the case that $k_4(k_3 \le x) \le [x]_4^3$, we have proved the following result.

<u>Theorem 5</u>. If $x = (\binom{n}{3} - 2(n-2))$, $[x]_4^3 - 1 \le k_4(k_3 \le x) \le [x]_4^3$, n > 6.

We believe that these Turan graphs, T(n, n - 2) are the maximal examples for n > 6; thus we make the following conjecture.

<u>Conjecture</u>. If n > 6, and $x = \binom{n}{3} - 2(n-2)$, then $k_4(k_3 \le x) = [x]_4^3 - 1$.

It is to be stressed that the above results and conjecture only apply to $K_3$, $K_4$ subgraphs. Presumably the results would be different for other $K_r$, $K_s$ subgraphs with r < s. For example, for $K_3$, $K_5$ subgraphs, and the Turan graphs, T(n, n - 2) we have the following table that was generated in Mathematica using the methods in [ 3 ].

| n | $K_3$ subgraphs | $K_5$ subgraphs | K5K3 bound |
|---|---|---|---|
| 6 | 12 | 0 | 1 |
| 7 | 25 | 4 | 6 |
| 8 | 44 | 20 | 23 |
| 9 | 70 | 61 | 65 |
| 10 | 104 | 146 | 151 |
| 11 | 147 | 301 | 307 |
| 12 | 200 | 560 | 567 |
| 13 | 264 | 966 | 974 |
| 14 | 340 | 1572 | 1581 |
| 15 | 429 | 2442 | 2452 |

This suggests that the number of $K_5$ subgraphs in T(n, n - 2) given $x = \binom{n}{3} - 2(n-2)$, $K_3$ subgraphs is (n -

5) less than the Kruskal-Katona bound $[x]_5^3$. The number of $K_5$ subgraphs of $K_n$ is $\binom{n}{5} - 2\binom{n-2}{3} + (n-4)$; while the K-K bound $[x]_5^3 = \binom{n-1}{5} + \binom{n-4}{4} + \binom{n-5}{3}$. On subtracting $(n-5)$ from $[x]_5^3$ and simplifying, it must be shown that: $\binom{n}{5} - 2\binom{n-2}{3} + (2n-9) = \binom{n-1}{5} + \binom{n-4}{4} + \binom{n-5}{3}$. This can be verified, in general, using a computer algebra system, such as Mathematica. Hence we have proved the following theorem.

<u>Theorem 6</u>. If $x = \binom{n}{3} - 2(n-2)$, then $k_5(k_3 \leq x) \geq \binom{n}{5} - 2\binom{n-2}{3} + (n-4)$.

Whether or not $k_5(k_3 \leq x) = \binom{n}{5} - 2\binom{n-2}{3} + (n-4)$, for some or all values of $n > 6$, is an open problem. Other combinations also yield interesting results that can be investigated.

4. <u>Some remarks on finding the number of $K_s$ subgraphs of a graph</u>. If we want to find the number of $K_s$ subgraphs of a connected graph $G$, it can be advantageous to first prune the graph by recursively removing all vertices, together with their edges, of degree less than $(s-1)$. The graph, $G^{s-1}$, that results from this pruning is known as the $(s-1)$-*core* of $G$ and this pruning operation is linear in the number of vertices and edges (see [6]). Clearly, no $K_s$ subgraph will be removed in obtaining the $(s-1)$-core, since each vertex in $K_s$, along with all its neighbors, is of degree, $(s-1)$. $G^{s-1}$ may be a smaller graph than $G$ with fewer $K_r$ subgraphs for $r < s$. Also, this can lead to better Kruskal-Katona estimates for upper bounds for $K_s$.

<div align="center">References</div>